\documentclass[12pt]{amsart}
\usepackage[margin=1in]{geometry}
\usepackage{amsthm}
\usepackage{times,eucal,fullpage,times,amsmath,amsthm,amssymb,mathrsfs,stmaryrd,enumerate,accents}
\usepackage[all]{xy}
\usepackage[usenames,dvipsnames]{color}
\usepackage{hyperref} 
\usepackage{epigraph}
\usepackage{graphicx}
\usepackage{fancyhdr} 
\usepackage{enumitem}
\newcommand{\be}{\begin{equation}}
\newcommand{\ee}{\end{equation}}
\newcommand{\bes}{\begin{equation*}}
\newcommand{\ees}{\end{equation*}}
\newcommand{\bea}{\begin{eqnarray}}
\newcommand{\eea}{\end{eqnarray}}
\newcommand{\beas}{\begin{eqnarray}}
\newcommand{\eeas}{\end{eqnarray}}
\newcommand{\ben}{\begin{note}}
\newcommand{\een}{\end{note}}
\newcommand{\bexl}{\vskip0.1em\noindent\hrulefill\vskip1em\begin{ExerciseList}}
\newcommand{\eexl}{\end{ExerciseList}\hrulefill}

\newcommand{\bthm}{\begin{theorem}}
\newcommand{\ethm}{\end{theorem}}
\newcommand{\bpro}{\begin{prop}}
\newcommand{\epro}{\end{prop}}
\newcommand{\bcor}{\begin{corollary}}
\newcommand{\ecor}{\end{corollary}}
\newcommand{\bcon}{\begin{conjecture}}
\newcommand{\econ}{\end{conjecture}}
\newcommand{\bp}{\begin{proof}}
\newcommand{\ep}{\end{proof}}
\newcommand{\blem}{\begin{lemma}}
\newcommand{\elem}{\end{lemma}}
\newcommand{\bn}{\begin{note}}
\newcommand{\en}{\end{note}}
\newcommand{\benum}{\begin{enumerate}}
\newcommand{\eenum}{\end{enumerate}}
\newcommand{\bed}{\begin{defn}}
\newcommand{\eed}{\end{defn}}
\newcommand{\brem}{\begin{remark}}
\newcommand{\erem}{\end{remark}}

\newcommand{\btik}{\begin{tikzpicture}\begin{axis}[scale=0.5,axis y line=center, axis x line=middle]}
\newcommand{\etik}{\end{axis}\end{tikzpicture}}

\let\into=\hookrightarrow

\let\cong=\equiv

\newcommand{\upperRomannumeral}[1]{\uppercase\expandafter{\romannumeral#1}}

\newtheorem{theorem}[equation]{Theorem}      
\newtheorem{lemma}[equation]{Lemma}          %
\newtheorem{corollary}[equation]{Corollary}  
\newtheorem{proposition}[equation]{Proposition}

\theoremstyle{definition}
\newtheorem{conj}[equation]{Conjecture}
\newtheorem*{example}{Example}

\theoremstyle{definition}
\newtheorem{defn}[equation]{Definition}
\theoremstyle{remark}

\theoremstyle{definition}
\newtheorem{remark}[equation]{Remark}

\numberwithin{equation}{section}

\let\into=\hookrightarrow
\let\isom=\simeq

\let\tensor=\otimes

\newcommand{\F}{{\mathbb F}}

\newcommand{\gal}{{\rm Gal}}

\newcommand{\Q}{{\mathbb Q}}

\newcommand{\Z}{{\mathbb Z}}

\renewcommand{\int}{\operatorname{int}}
\renewcommand{\O}{{\mathcal O}}

\renewcommand{\wp}{{\mathfrak p}}

\newcommand{\fq}{{\mathbb{F}_q}}
\newcommand{\fqt}{\fq(t)}
\newcommand{\bQ}{\bar{\Q}}
\renewcommand{\bpro}{\begin{proposition}}
\renewcommand{\epro}{\end{proposition}}

\long\def\comment#1\endcomment{}

\begin{document}

\title[]{Methods for constructing elliptic and hyperelliptic curves with rational points}%
\author{Kirti Joshi}%
\address{Math. department, University of Arizona, 617 N Santa Rita, Tucson
85721-0089, USA.} \email{kirti@math.arizona.edu}

\thanks{}%
\subjclass{}%
\keywords{}%


\begin{abstract}
I provide methods of constructing elliptic and  hyperelliptic curves over global fields with  interesting rational points over the given fields or over large field extensions. I also provide a elliptic curves defined over any given number field equipped with a rational point, (resp. with two rational points) of infinite order over the given number field,  and elliptic curves over the rationals with two rational points over `simplest cubic fields.' I also provide hyperelliptic curves of genus exceeding any given number over any given number fields with points (over the given number field)  which span a subgroup of rank at least $g$ in the group of rational points of the Jacobian of this curve. I also provide a method of constructing hyperelliptic curves over rational function fields with rational points defined over field extensions with large finite simple Galois groups, such as  the Mathieu group $M_{24}$. 
\end{abstract}
\maketitle
\epigraph{My scheme is far more subtle.\\ Let me outline it for you.}{Bertie Wooster in \cite{wodehouse}}

\tableofcontents

\section{Introduction}
This note and \cite{joshi17-legendre}, grew out of my attempts to understand Ulmer's remarkable construction of rational points on the Legendre elliptic curve over $\fq(t)$ after his insightful colloquium talk at the University of Arizona (on March 1,2017). While \cite{joshi17-legendre} was mostly about constructing rational points on Legendre elliptic curves, this paper describes other methods of constructing elliptic and hyperelliptic curves with rational points.

The beginning of this construction method dates back to around 1989. I discovered this construction for $K=\Q$  and genus one in 1989. In 1988 Don Zagier taught a course on Gross-Zagier Theorem at the Tata Institute of Fundamental Research and around that time, and in response to  a question posed by C.~S.~Dalawat during one coffee-table conversation, I wrote down the following elliptic  curve
\be\label{eq:orig} 
E: y^2=x^3+x+d^2 
\ee
\renewcommand{\thefootnote}{\fnsymbol{footnote}}
which I showed for any integer $d\geq 3$  (note that the discriminant $4+27d^4$ is obviously never zero) has rank at least one by showing that the manifest point $(0,d)$ has infinite order (for any $d\geq 3$). This is immediate from the Lutz-Nagell Theorem: the point is not of order two as $d\neq 0$ and so it has order at least three and in which case $d^2$ must divide $4+27d^4$ which it cannot  as this implies $d^2|4$, so $d\leq2$ which contradicts our hypothesis on $d$. 
Moreover by a simple search, this curve also has rank at least one for $d=\pm1,\pm2$ and hence it has rank at least one for all $0\neq d\in\Z$. 
Two curves associated with $d,d'\in\Q^*$,  are  $\bQ$-isomorphic if and only if $d/d'$ is a fourth root of unity. In particular this is not a family of twists of any elliptic curve over $\Q$.

The following table of values shows the ranks of these curves for  $3\leq d\leq 100$\footnote{This data has been computed recently using \cite{sage}, my own computations, for a more modest range of  values ($d\leq 20$), in 1989--92 were carried out with a patchwork of small search programs, with many stretches of manual and symbolic computations as well, were written mostly in PASCAL for VAX and CYBER ``mainframes'' which were owned by the Tata Institute at the time. When John Cremona's `mwrank' became available I was able to extend these computations.}. Apart from curves of rank one and two, there are curves of rank three and one of rank four for $d=82$ which is striking as the coefficients of this curve $y^2=x^3+x+82^2=x^3+x+6724$ are not too big for a curve of this rank.
\begin{center}
	\begin{tabular}{|c|l|}
		\hline rank & $d$  \\ 
		\hline 1  & $5, 12, 17, 32, 35, 36, 39, 40, 41, 42, 48, 49, 52, 54, 55,\cdots$ \\ 
		\hline 2  & $3, 6, 7, 8, 9, 10, 11, 13, 15, 16, 19, 20, 22, 24, 25, 26, 27, 28, 29, 30,\cdots$ \\ 
		\hline 3  &  $14, 18, 21, 23, 31, 44, 57, 68, 69, 74, 78, 79,\cdots$ \\ 
		\hline 4  & $82,\cdots$ \\
		\hline
	\end{tabular} 
\end{center}
For $d=82$ the four generators are of the Mordell-Weil group are $(0 : 82 : 1), (12 : 92 : 1), (60 : 472 : 1), (465/4 : 10049/8 : 1)$.

A few days after I discovered (\ref{eq:orig}), I found the curve
\be\label{eq:orig2} 
E: y^2=(x-a)(x-b)(x-c)+d^2
\ee
which has six manifest points $(a,\pm d),(b,\pm d),(c,\pm d)$, and   for a suitable set of choices of $a,b,c,d$, the points  $(a,d),(b,d)$ have infinite order and are linearly independent in the Mordell-Weil group  (note that the three points $(a,d),(b,d),(c,d)$ are on the line $y=d$ and so there is one relation between these three points (a variant of \eqref{eq:orig2} with rank at least two) is described in my paper with Pavlos Tzermias \cite{joshi99}). I did not publish my results but I did share them with some of my colleagues (notably with C.~S.~Dalawat who had raised the question) and also described them in my correspondence with C.~B.~Khare but I did not think about these issues for the next few years, except to observe, when Kolyvagin's work became available, that my examples also provided families of $L$-functions with analytic ranks $\geq 1$ and $\geq 2$ respectively.

After I writing the present note I found that in recent elliptic curve literature, special cases of \eqref{eq:orig2} have been studied. For example taking $a=0,b=m,c=-m$ gives the curve (of rank at least two)
$$y^2=x^3-m^2x+d^2,$$
studied by \cite{fujita17}, while in \cite{izadi17} this special case is studied with a view to rank three sub-families. Let me also say that by mentioning above history of \eqref{eq:orig} and \eqref{eq:orig2} it is not my intent to claim any priority over these  authors (or others who may have found this equation too);  they, too,   discovered   independently that these curves have interesting properties.

Here is brief description of the new results which grew out of my attempts to understand Doug Ulmer's remarkable construction of explicit rational points in characteristic $p>0$: in Theorem~\ref{th:orig-2} I show that a variant of \eqref{eq:orig} also provides families of elliptic curves over any number field with a rational point of infinite order defined over that field. In Theorem~\ref{th:orig-2} I show that a variant of \eqref{eq:orig2} provides two rational points, both of infinite order, defined over any given number field. These examples are  well-adapted for studying the conjecture of \cite{denef80} thanks to a result of \cite{poonen02} (see Conjecture~\ref{con:denef}). 

In Theorem~\ref{th:joshi-tzermias-2}, which is a variant of the result of \cite{joshi99}, I provide examples of infinite families of hyperelliptic curves defined over any given number field with many rational points and Jacobian rank at least as large as the genus.  

In Theorem~\ref{th:main0} I provide a general method for constructing hyperelliptic curves with rational points over large fields. This has a number of consequences. For instance in Theorem~\ref{th:symmetric} I  construct hyperelliptic curves over number fields with about $4g$ rational points over finite extensions of large degree (about $g^g$)--it is in this context that I stumbled upon the result of \cite{osada87} (for other uses of this result see \cite{joshi17-legendre}); a variant of this  constructs hyperelliptic curves over $\Q$ with rational points over a finite extension with Galois group $A_n$ for all odd $n$ (see Theorem~\ref{th:An}); and curves over $\fq(t)$ with points over a finite extension of $\fq(t)$ with Galois group equal to ${\rm PGL}(n,\fq)$ (see Theorem~\ref{th:pgl}); and also curves over  $\F_2(t)$ with rational points over a finite extension of $\F_2(t)$ with Galois group equal to Mathieu group $M_{24}$ (see Theorem~\ref{th:m24}). There are other examples possible but I believe this list is already quite interesting; except for the Mathieu Group case, in all other cases one has an infinite family of curves of these types.

I do not claim that my examples break any records and I certainly do not claim that my methods or examples  are optimal in any way but I do hope to convince the reader that the examples which I describe here are remarkable for their simplicity, and versatility and that they will be useful to any student of rational points on curves and can be the starting point for other investigations. For example, one natural problem which suggests itself is that these higher genus curves could be used to test Minhyong Kim's conjecture that there is some higher dimensional Selmer variety in which the image of these points spans a proper subspace. Number of other conjectures are also included here which are based on my computations.

Conjecture~\ref{con:main} is a somewhat optimistic conjecture which predicts that the rational points constructed in the hyperelliptic curve examples of Theorem~\ref{th:main0} have infinite order. Some numerical evidence for this is provided in \ref{ss:examples} and Conjecture~\ref{con:main} is proved in some cases in Theorem~\ref{th:main1}.

In Subsection~\ref{ss:cubic} I return to elliptic curves again and provide examples of elliptic curves over $\Q$ with two rational points over ``simplest cubic fields'' (see \cite{shanks74} for this terminology). Using methods of Theorem~\ref{th:main1}  I show that these two points are linearly independent in the Mordell-Weil group in many cases (see Theorem~\ref{th:main2}). In Conjecture~\ref{con:orthogonality} I state a general conjecture about these curves which I found experimentally.   I provide a number of numerical examples illustrating this conjecture.

It is a pleasure to thank Douglas Ulmer for conversations and correspondence and for his colloquium talk which led to this paper (and \cite{joshi17-legendre}). I would also like to take this opportunity to thank C.~S.~Dalawat for asking  questions which led me to construct the curves \eqref{eq:orig} and \eqref{eq:orig2} a long time ago. It is a pleasure to thank Tanmoy Bhattacharya for many hours of tutorials on programming on VAX and CYBER he provided me in 1987--88 which made my computations on \eqref{eq:orig} possible. Thanks are also due to Kapil Paranjape for many conversations about elliptic curves and programming on VAX/CYBER during that period.

\numberwithin{equation}{subsection}

\section{Elliptic Curves}
\renewcommand{\wp}{\mathfrak{p}}
\newcommand{\wq}{\mathfrak{q}}
\newcommand{\ok}{\O_K}
\subsection{Elliptic curves with a point of infinite order over a given number field}
In this section, I provide a simple construction of an elliptic curve over every number field with Mordell-Weil group of rank at least one over this number field.  Let $K$ be any number field, $\ok$ be the ring of integers of $K$. For prime ideals $\wp,\wq\neq0$ of $\ok$ let $p$ (resp. $q$) be the corresponding primes of $\Z$ lying below $\wp$ (resp. $\wq$). 

\bthm\label{th:orig-2}
Let $K$ be any number field as above and $\wp,\wq$ be non-zero prime ideals of $\ok$ not lying over $(2)\subset \Z$, with residue fields of prime order $p$ (resp. $q$) and $q>p+1+2\sqrt{p}$. Let $0\neq d\in \ok$ be chosen so that
\benum[label={\bf(\arabic{*})}]
\item\label{th:orig-2:hyp-1} $K=\Q(d)$,
\item\label{th:orig-2:hyp-2} $16(4+27d^4)\not\in\wp$,
\item\label{th:orig-2:hyp-3} $d\in\wq$.
\eenum
Then 
\benum[label={\bf(\arabic{*})}]
\item $P=(0,d)$ is a point of infinite order on the elliptic curve $$E:y^2=x^3+x+d^2,$$
\item  neither $E$ nor $P$ can be defined over any proper subfield contained in $K$. 
\item for $d,d'\in\ok$, the corresponding curves are isomorphic over $\bQ$ if and only if $d=d'\zeta_4$ for some fourth root of unity $\zeta_4$. 
\item In particular if $d\neq d'\zeta_4$ for any fourth root of unity, then the corresponding curves are not twists of each other.
\eenum
\ethm 

\bp 
First let me show that the hypothesis of the theorem are satisfied infinitely often. By the proof of the primitive element theorem one sees that there exits (infinitely many) $d_1\in\ok$ such $K=\Q(d_1)$. So suppose one has chosen a $d_1\in\ok$ such that $K=\Q(d_1)$. Then note that any rational integer multiple of $d_1$ has the same property. It is well-known that for any number field there is a set of primes ideals of density one  with prime residue fields (see By Lemma~\ref{le:density-lemma}). So there are infinitely many prime ideals $\wp,\wq$ which satisfy the residue field hypothesis of the theorem.  Choose  primes $\wp,\wq$ lying over odd rational primes $p,q$ such that $q>p+1+2\sqrt{p}$ and $\ok/\wp=\F_p,\ok/\wq=\fq$. Choose an integer $0\neq a\in\Z$ such that $a\in\wp\cap\wq=\wp\wq$. Let $d=ad_1$. Then $K=\Q(d)$, $d\in\wq$ and $-16(4+27d^4)\not\in\wp$ (otherwise $4\in\wp$). Thus there is an infinite set of $(d,\wp,\wq)$ which satisfy all the hypothesis of the theorem.  Note further that if $16(4+27d^4)\not\in\wp$ then $16(4+27d^4)\neq0$ so $E$ is always an elliptic curve! Note that $E$ has good reduction at $\wp$ as the discriminant is not contained in $\wp$ by hypothesis and as $d\in\wq$ but $\wq\not|2$ so $\wq$ also does not divide the discriminant $16(4+27d^4)$ of $E$. In particular one sees, as $d\in\ok$,  that the equation for $E$ is also minimal at $\wp$ and $\wq$.

Suppose $P=(0,d)$ is a point of finite order $m>1$. Then $m\geq 3$: for if $2P=0$ then by \cite[Lutz-Nagell Theorem]{silverman-arithmetic} one has $d=0$ which is not the case by my assumptions on $d$. If $(q,m)=1$ then reducing modulo $\wq$ one sees (by \cite[Chap 7, Proposition 3.1]{silverman-arithmetic}) that $E[m]\into E(\ok/\wq)$ and image of $P$ which is $(0,0)\bmod{\wq}$ is of order two on the reduction $y^2=x^3+x\bmod{\wq}$. So again $m=2$ or $(q,m)=q>1$. As $m\neq 2$, let $m=q^an$ with $(q,n)=1$ then $Q=nP$ is a point of order $q^a\geq q$. Now reducing modulo $\wp$ and noting that $E[q^a]\into E(\ok/\wp)$ one sees that 
$$q\leq|E[q^a]|\leq |E(\ok/\wp)|\leq p+1+2\sqrt{p}$$
 which is a contradiction as $q>p+1+2\sqrt{p}$. So $P$ has infinite order in $E(K)$.
 
Now suppose $d,d'\in K$ then the corresponding  curves are isomorphic over $\bQ$ if and only if they have the same $j$-invariants. Using the formula for $j$-invariants one sees that this means $$\frac{4}{16(4+27d^4)}=\frac{4}{16(4+27d^{'4})}$$ from which the last claim follows. In particular one sees that curves for $d$ and $d'$ are not twists of each other if $d\neq d'\zeta_4$ for any fourth-root of unity $\zeta_4$.
\ep

The following lemma is certainly well-known (see \cite{lang-algebraic}) and ensures that for any number field $K$ there are infinitely many $\wp,\wq$ satisfying the hypothesis of the theorem.

\blem\label{le:density-lemma}
Let $K/\Q$ be a number field. Then the set of primes $\wp\subset\ok$ of $\ok$  such that $\ok/\wp=\F_p$ is of Dirichlet density one in $K$.
\elem

\subsection{A numerical example}
Here is a numerical example illustrating the constructions of Theorem~\ref{th:orig-2}. 
\newcommand{\exd}{\sqrt{-647}} 
Let $K=\Q(\exd)$ then class number of $K$ is $h_K=23$. Let $\wp=\left(29,\frac{7+\exd}{2}\right)$ and $\wq=\left(67,\frac{31+\exd}{2}\right)$ both are prime ideals. Then clearly $67>29+1+2\sqrt{29}$. Choose $d=98+\exd=67+2\left(\frac{31+\exd}{2}\right)\in\wq$ and one checks that $-16(4+27d^4)\not\in\wp$. Then according to the theorem 
$$E:y^2=x^3+x+d^2=x^3+x+\left(98+\exd\right)^2$$
has a rational point $P=(0,d)=\left(0,98+\exd\right)$ of infinite order. One finds using \cite{sage} that torsion of $E(K)$ is zero and the N\'eron-Tate height of $P$ is $\hat{h}(P)\approx3.54047$. Thus $P$ is non-torsion (as predicted by the theorem).

\subsection{Elliptic curves  with two points of infinite order over a given number field}
Let me also mention the following variant of Theorem~\ref{th:orig-2} which will be used in the following subsection. 
\bthm\label{th:orig-3}
Let $K$ be any number field as above and $\wp,\wq$ be non-zero prime ideals of $\ok$ not lying over $(2)\subset \Z$, with residue fields of order $p$ (resp. $q$) and $q>p+1+2\sqrt{p}$. Let $0\neq d,\beta\in \ok$ be chosen so that
\benum[label={\bf(\arabic{*})}]
\item\label{th:orig-3:hyp-1} $K=\Q(d)=\Q(\beta)$,
\item\label{th:orig-3:hyp-2} $16(-4\beta^6+27d^4)\not\in\wp$,
\item\label{th:orig-3:hyp-3} $d\in\wq$,
\item\label{th:orig-3:hyp-4} $\beta\not\in\wp$,
\item\label{th:orig-3:hyp-5} $\beta\not\in\wq$.
\eenum
Then 
\benum[label={\bf(\arabic{*})}]
\item $P=(0,d)$ is a point of infinite order on the elliptic curve $$E:y^2=x^3-\beta^2 x+d^2,$$
\item $Q=(\beta,d)$ is also a point of infinite order on $E$, and
\item  neither $E$ nor $P,Q$ can be defined over any proper subfield contained in $K$. 
\eenum
\ethm 
Note that the equation $E:y^2=x^3-\beta^2 x+d^2=x(x-\beta)(x+\beta)+d^2$ is special case of \eqref{eq:orig2} and as mentioned earlier, for $K=\Q$,  Mordell-Weil ranks of the  curves $y^2=x^3-\beta^2x+d^2$ have been recently considered by \cite{fujita17} and \cite{izadi17}; they provide sub-families of these curves with ranks three and four. The proof of Theorem~\ref{th:orig-2} is similar to the proof of Theorem~\ref{th:orig-3} and is left to the reader.

\subsection{Another numerical example}
\renewcommand{\exd}{\sqrt{94546}}
Here is a numerical example illustrating the constructions of Theorem~\ref{th:orig-3}. 
Let $K=\Q(\exd)$ then class number of $K$ is $h_K=80$. Let $\wp=\left(29,8+\exd\right)$ and $\wq=\left(67,3+\exd\right)$ both are prime ideals. Then clearly $67>29+1+2\sqrt{29}$. Choose $d=5905-265\exd\in\wq$ and $\beta=-104-195\exd$. Then  one checks that $-16(4(-\beta^2)^3+27d^4)\not\in\wp$  and that one has $\beta\not\in\wp\wq$. Then according to the theorem 
$$E:y^2=x^3-\beta^2x+d^2=x^3-(3595122466+40560\exd) x+\left(5905-265\exd\right)^2$$
has two rational points $$P=(0,d)=\left(0,5905-265\exd\right), Q=(\beta,d)=(-104-195\exd,5905-265\exd).$$  One finds using \cite{sage} that torsion of $E(K)$ is zero and the N\'eron-Tate height pairing for  $P,Q$ is 
\be
\begin{pmatrix}
 10.6548595059767  & -5.32742744032818\\
-5.32742744032818 & 10.9147767670134
\end{pmatrix}
\ee
whose determinant is $\approx87.9139298596650$. Thus $P,Q$ are both non-torsion and linearly independent (as predicted by the theorem).

\subsection{On Denef's conjecture}
In \cite{poonen02} it  was  shown that a conjecture of \cite{denef80} on Diophantine definability of $\Z\subset \O_L$ follows from a conjecture about elliptic curves of rank one. Specifically it was  shown that  for  number fields $K\subset L$,  if there exists an elliptic curve over $K$ which has rank exactly equal to one over $K$ and $L$ then there exists a Diophantine definition of $\ok$ over $\O_L$. In particular these works imply that Hilbert's Tenth problem is undecidable over $\O_L$ (if such an elliptic curve is found for suitable number fields $K$, see \cite{poonen02} and \cite{shlapentokh08}. The condition that the rank is one was later replaced by  the condition that there is an elliptic curve over $K$ such that the ranks over $L$ and $K$ coincide and this number is positive (see \cite{shlapentokh08}). In \cite{mazur10} this conjectured existence of rank preserving elliptic curves of positive rank is proved under a certain parity conjecture (also see  \cite{murty15} for related results based on analytic considerations notably for totally real fields).

Theorem~\ref{th:orig-2} and Theorem~\ref{th:orig-3}, which produce a rich supply of elliptic curves over any given number field with an explicit points of infinite order, reduce the question of finding rank preserving elliptic curves of positive rank (and hence  Hilbert's Tenth Problem for $\O_L$)  to the following more concrete assertions (see Conjectures~\ref{con:denef}, \ref{con:denef2} and \ref{con:denef3}) which may perhaps be accessible. At any rate proof of these conjectures, for a given $K$ and $L$, implies that $\ok$ has a Diophantine definition in $\O_L$.  In particular under the hypothesis of Theorem~\ref{th:orig-2} one has positivity of the rank (for the elliptic curves I construct). So rank positivity hypothesis in \cite{shlapentokh08} holds for these curves. Since Theorem~\ref{th:orig-2} provides an abundance of curves over a given number field $K$ with positive rank and an explicit point, it seems reasonable to expect that for any finite cyclic extension $L/K$ there is at least one (and possibly infinitely many) of these curves which preserve their ranks over $L$.  I provide here two explicit versions of this conjecture--the second may be more difficult to prove (but for more on this see \cite{mazur10}):
\begin{conj}\label{con:denef}
Let $K\subset L$ be a cyclic extension of prime degree of number fields and suppose $\sqrt{-1}\not\in L$. Then there exists a $0\neq d\in\ok$ satisfying the hypothesis \ref{th:orig-2:hyp-1}--\ref{th:orig-2:hyp-3} of Theorem~\ref{th:orig-2} such that the point $(0,d)$ on $E:y^2=x^3+x+d^2$ 
is of infinite order and generates $E(L)\tensor\Q$.
\end{conj}

\begin{remark} 
Let me remark that the restriction to cyclic extensions of prime degree in the above conjecture comes from the observation of \cite{mazur10}.
Secondly if $i=\sqrt{-1}\in L$ then there is another rational point $(i,d)$ and this is often linearly independent of $(0,d)$. For example for $y^2=x^3+x+49$, the points $P=(0,7)$ and $Q=(i,7)$ are linearly independent.
\end{remark}

\begin{conj}\label{con:denef2}
Let $K\subset L$ be a pair of number fields. Then there exists a $0\neq d\in\ok$ satisfying the hypothesis \ref{th:orig-2:hyp-1}--\ref{th:orig-2:hyp-3} of Theorem~\ref{th:orig-2} such that ${\rm rk}(E(L))={\rm rk}(E(K))\geq 1$.
\end{conj}

To avoid the hypothesis $\sqrt{-1}\not\in L$ consider more generally the curve $y^2=x^3+\alpha x+d^2$. This has a rational point $P=(0,d)$ and the proof of Theorem~\ref{th:orig-2} also shows that this point is of infinite order (assuming for instance that $\alpha\not\in\wq$). This curve also has another manifest point: $(\sqrt{-\alpha},d)$ and this point is also of infinite order in $E(K)$ if and only if $\sqrt{-\alpha}\in K$. Given number fields $K\subset L$ it is certainly possible to choose $\alpha\in K$ such that $\sqrt{-\alpha}$ is not in $L$. On the other hand by choosing $\alpha=-\beta^2$ suitably  one sees from Theorem~\ref{th:orig-3} that $(\beta,d)$ is a $K$-rational point of infinite order.  So the more general version of the conjecture is the following conjecture.

\begin{conj}\label{con:denef3}
Let $K\subset L$ be an extension of number fields. Then there exists an $\alpha\in\ok$ and $0\neq d$ in $\ok$ satisfying the hypothesis \ref{th:orig-2:hyp-1}--\ref{th:orig-2:hyp-3} of Theorem~\ref{th:orig-2} such that
$$1\leq {\rm rk}(E(K))={\rm rk}(E(L))\leq 2.$$
\end{conj}


\section{Hyperelliptic Curves}\label{se:high-genus}
\subsection{Hyperelliptic Curves over any number field with Jacobian of rank at least $g$}
In 1999 during the Arizona Winter School, where Coleman-Chabauty methods were being applied to specific genus two curves, I showed my example (\eqref{eq:orig}) to Pavlos Tzermias, (both of us were post-doctoral faculty at Arizona at the time) and suggested that we could try and study the genus two case and see if we could prove that the point $(0,d)$ has infinite order in the Jacobian of that curve. This led to our paper \cite{joshi99} where we proved substantially more than either of us had initially hoped to prove. We showed that for any odd prime $p\geq 5$ and any integers $a_1,\ldots,a_p$ representing the $p$ distinct residue classes modulo $p$ the curve:
\be\label{eq:joshi-tzermias} 
y^2=(x-a_1)\cdots (x-a_p)+p^2d^2
\ee
has at least $2p$ or about $4g$ rational points $(a_i,\pm pd)$ and its Jacobian has Mordell-Weil  rank at least $g$ over $\Q$. This example, despite its simplicity and generality, has largely been forgotten in the Coleman-Chabauty method and rational points literature (at least as far as I can tell). Note that in \cite{joshi99} we did not use Coleman-Chabauty method to compute ranks--rather we turned the method on its head and proved that ranks must be ``large'' as we have more rational points than allowed by Coleman-Chabauty bound. Our paper also has a genus one variant of \eqref{eq:orig2} with rank at least two.   It is amusing to note that  for $d=1$ in \eqref{eq:orig} one gets the curve $y^2=x^3+x+1$ which also has rank one (though this needs more work than Lutz-Nagell Theorem argument which I sketched above) this genus one curve is the quotient of the genus two curve $y^2=x^6+x^2+1$, with the solution $(1/2,9/8)$, which occurs in the work of Diophantus and one finds a complete list of rational points on  it in J.~L.~Wetherell's  thesis  (see \cite{wetherell99}). During AWS 1999, Wetherell explained his approach to rational points on the Diophantus curve using Coleman-Chabauty method and my work with Tzermias was to some extent prompted by Wetherell's lectures.

Here is a variant of \cite{joshi99} which provides a family of curves of genus $g=(p-1)/2\geq 2$ for a positive density of odd primes $p\geq 5$, over any number field, and with a collection of about $4g$ points which generate a subgroup of rank at least $g$. 

\bthm\label{th:joshi-tzermias-2}
Let $K/\Q$ be any number field. Let $\wp$ be a  prime of $\ok$ lying over an unramified rational prime $p\geq 5$. Let $\ok/\wp=\F_q$, so $q=p^m$ for some $m\geq 1$. Let $a_1,\ldots,a_q\in\ok$ belong to the $q$ distinct residue classes modulo $\wp$. Let $0\neq d\in\wp$ and suppose that $K=\Q(a_1,\ldots,a_p,d)$. Consider the hyperelliptic curve
$$C:y^2=(x-a_1)\ldots (x-a_q)+d^2.$$
Then
\benum[label={\bf(\arabic{*})}]
\item  $C$ has $2q+1$ $K$-rational points given by $(a_i,\pm d)$ and $\infty$,
\item  its Jacobian $J$ has $K$-rank at least $g$. 
\item Moreover this curve cannot be defined over a subfield of $K$. 
\item In particular for a set of primes $p\geq 5$ of positive density, there is a family of curves of genus $g=\frac{p-1}{2}$ which is defined over $K$ and such that the Mordell-Weil rank is at least $g$.
\eenum
\ethm
\bp 
The proof is same as the one given in \cite{joshi99} I recall it here for convenience: suppose if possible that the rank of $J(K)$ is at most $g-1$. By construction there are at least $2q+1$ $K$-rational points given by $\infty$ and $(a_i,\pm d)$. The latter reduce to $(\alpha,0)$ for every $\alpha\in\ok/\wp=\F_q$. So modulo $\wp$ the curve, which is now given by $y^2=x^q-x$, has exactly $q+1$ rational points over $\F_q=\ok/\wp$.   On the other hand,  as the genus of the curve is $\frac{q-1}{2}$, and by our supposition the rank is at most $g-1$, so by \cite{coleman85} the number of $K$-rational points  $$2q+1\leq \#C(K)\leq2g-2+\#C(\ok/\wp)=q-1-2+q+1=2q-2.$$ So one has arrived at a contradiction. Hence the  rank of $J(K)$ is at least $g$.

The assumption that the field generated by $a_1,\ldots,a_q,d$ is $K$ ensures that the points are not all defined over a subfield of $K$. The last assertion follows from Lemma~\ref{le:density-lemma} upon considering primes $\wp$ such that $\ok/\wp=\F_p$.
\ep

\brem
In \cite{joshi99} we showed that for $K=\Q$ if the polynomial $(x-a_1)\cdots (x-a_p)+d^2$ is irreducible then the rank is at least $2g$. It seems reasonable that some version of this result  of \cite{joshi99} also holds for any number field $K$. In fact ranks  can be even higher than $2g$ (in this situation). For instance in \cite{gallegos} it is shown that the curve ($K=\Q$, $p=5$ and $g=2$ for this example)
$$y^2=(x + 19)(x + 20)(x + 21)(x + 22)(x + 23) + 20^2$$
has Jacobian of $\Q$-rank at least five.
\erem

\subsection{Hyperelliptic curves  with points over fields of large degree}
I describe a  method which provides an interesting collection of rational points on higher genus curves. Let $K$ be a field of characteristic $p\neq 2$ ($p=0$ is allowed). If $C/K$ is a smooth, projective curve, let $J_C=Jac(C)$ be the Jacobian of $C$ and if $C(K)\neq \emptyset$, then embed $C\into J$ by in the usual way one of the $K$-rational points. I construct curves over small global fields with points over rather large fields (typically the examples I construct have genus $g$ and have points at least $4g$ rational points over fields of degree $\approx g^g$ and these points are not defined over any proper sub-extension which is Galois over the base-field). My constructions are based on the following simple idea and provides curves over $\Q$ with rational points over fields of very large degrees in comparison with their genus. 

For a finite group $G$, let $m(G)$ be the minimum of dimensions of $\Q$-rational, faithful representations of $G$ on a $\Q$-vector space.
\bthm\label{th:main0}
Let $f(X)\in K[X]$ be an irreducible polynomial of odd degree  with splitting field $L_f/K$. Suppose that for  some  $d\in K^*$ the polynomial 
$$g(X)=
f(X)+d^2 
$$
has distinct roots. 
Then 
\benum[label={\bf(\arabic{*})}] 
\item\label{th-points-1} Let $f(u)=0$ be a root of $f(X)$. Then $(u,\pm d)$ is an $L$-rational point on $$C:y^2=g(X).$$
\item\label{th-points-2} In particular $|C(L)|\geq 2\deg(f)\geq 4g+2$.
\item\label{th-points-3} If $\gal(L_f/K)$ is a simple group then $L_f$ is the smallest Galois extension of $K$ over which $(u,\pm d)$ is defined.
\item\label{th-points-4-trans} Moreover $\gal(L_f/K)$ acts transitively on the $\deg(f)$ points $P=(u,d)$ where $u$ runs through the roots of $f(X)=0$.
\item\label{th-points-5-rank} Suppose for one root $u$, $P=(u,d)$ is non-torsion in $J(L_f)$, then rank of the subgroup of $J_C(L_f)$ generated by these points
is at least  $m(\gal(L_f/K))$.
\eenum 
\ethm
\bp 
Let $u$ be a root of $f(X)$, let $G=\gal(L_f/K)$. Then clearly $P_\pm=(u,\pm d)$ is in $C(L_f)$. Genus of $C$ satisfies $\deg(f)=2g+1$. So $|C(L_f)|\geq 2\deg(f)+1=4g+2$. If $G$ is simple then there are no nontrivial quotients of $G$. So $L_f$ is the smallest Galois extension of $K$ over which $(u,\pm d)$ are defined. Since $f$ is irreducible, $G$ acts transitively on the roots of $f(X)$ in $L_f$ and this induces a transitive action of $G$ on the $\deg(f)$ points $(u,d)$ where $u$ runs through the roots of $f(X)$. If one $(u,d)$ is torsion then by the transitivity of Galois action all of them are torsion. So if one of them is non-torsion the so are all of them. In particular $G$ acts on a non-zero $\Q$-vector space. Further as $G$ is simple this action is faithful (the trivial action is not transitive). Hence the subspace generated by these points is at least as big as $m(G)$.

The transitivity of the action of $G=\gal(L_f/K)$ follows from the fact that $G$ action on the roots of $f(X)$ is transitive by the irreducibility of $f(X)$.

Now suppose $P=(u,d)$ is non-torsion in $J_C(L_f)$ then by the Galois action this is true for all the roots $u$ of $f(X)$. So $G$ acts on the non-zero vector space generated by these points. Thus the span of these points in $J_C(L)\tensor\Q$ is of dimension at least $m(G)$.
\ep

\subsection{Hyperelliptic Curves over $\Q$ with points over large extensions}
Here are some consequences of this simple idea.

In \cite{osada87} it was shown that there exists infinitely many trinomials $X^n+aX^s+b\in\Q[X]$ whose Galois group is $S_n$. An explicit example of such a polynomial was also given in \cite{osada87}: $X^n-X-1$ has Galois group $S_n$. 
\bthm\label{th:symmetric}
Let $f(X)=X^n+aX^s+b\in\Q[x]$ for $n\geq 5$. Let $L/\Q$ be the splitting field of $f(X)$ and assume that $\gal(L/\Q)=S_n$. Let $g(X)=f(X)+d^2$. Then there exist infinitely many $d\in \Z$ is  such that $g(X)$ has distinct roots.  Then
\benum[label={\bf(\arabic{*})}]
\item $\gal(L/\Q)=S_n$ and
\item $|C(L)|\geq 2n$ and hence
\item $|C(L)|>\frac{\log([L:K])}{\log\log([L:K])}$
\eenum
\ethm

\bp 
The discriminant of $g$ is a polynomial in $d$. So the existence of infinitely many $d$ is clear. Hence $C$ is a nonsingular for such $d$. Now suppose $n$ is sufficiently large then one can use Stirling's estimate for $n!=[L:K]$ to get the last assertion. 
\ep

In \cite{hermez01} one finds the construction of explicit polynomials of the form $X^n+aX^m+b\in\Q[X]$ whose Galois group over $\Q$ is the Alternating group $A_n$. So one may use this to produce hyperelliptic curves with primitive points.
\bthm\label{th:An}
Let $f(X)=X^n+aX^m+b\in\Q[X]$ be a polynomial whose splitting field $L_f/\Q$ has $\gal(L_f/\Q)=A_n$ and $n\geq 7$. Suppose $d\in\Q$ is chosen so that $g(X)=f(X)+d^2$ has distinct roots. Then $C:Y^2=g(X)$ has an $L_f$-primitive point $P=(u,\pm d)\in C(L)$. Moreover $|C(L)|\geq 4g+2$.
\ethm

\subsection{Hyperelliptic curves over function fields}
The method I have described here is also applicable to curves over function fields. Here is a very small sample of the sort of results one can prove. In \cite{abhyankar97} Abhyankar constructed trinomial equations over $K=\F_q(t)$ with large Galois groups. One can use Theorem~\ref{th:main0} to these polynomials and obtain hyperelliptic curves with many points over highly non-abelian extensions of $\fqt$
\bthm\label{th:pgl}
Let $p$ be an odd prime, and $q$ be a power of $p$. Let $m\geq 2$ be an integer. Let $\langle m-1\rangle=\frac{q^m-1}{q-1}$. Let $K=\fq(t)$ and let $L_f/K$ be the splitting field of $U^{\langle m-1\rangle}+U+t\in K[U]$. Suppose $d\in K^*$ be such that $g(X)=f(X)+d^2$ has distinct roots and let  $C:Y^2=g(X)$. Then
\benum[label={\bf(\arabic{*})}]
\item $\gal(L_f/K)=PGL(m,\fq)$,
\item $|C(L)|\geq 2\frac{q^m-1}{q-1}$.
\eenum
\ethm

For $p=2$ one obtains, using \cite{conway10}, the following
\bthm\label{th:m24}
Let $p$ be an odd prime, and $q$ be a power of $p$. Let $m\geq 2$ be an integer. Let $\langle m-1\rangle=\frac{q^m-1}{q-1}$. Let $K=\fq(t)$ and let $L_f/K$ be the splitting field of $U^{\langle m-1\rangle}+U+t\in K[U]$. Suppose $\alpha,d\in K^*$ be such that $g(X)=f(X)+d(d+\alpha)$ has distinct roots and let  $C:Y^2+\alpha Y=g(X)$. Then
\benum[label={\bf(\arabic{*})}]
\item $\gal(L_f/K)=M_{24}$ (Mathieu Group),
\item $|C(L)|\geq 48$,
\item $L$ is the smallest Galois extension of $K$ which contains $(u,d)$ where $u$ is any root of $U^{\langle m-1\rangle}+U+t=0$.
\eenum
\ethm

\subsection{The rank conjecture}
Let $C/K$ be a  one of the curves described in the preceding sections (in particular $C$ is hyperelliptic, equipped with a $K$-rational point). Let $J_C=Jac(C)$ be the Jacobian of this curve and embed $C\into J_C$ using a $K$-rational point $P_0=\infty$. Now suppose $L/K$ is a finite Galois extension and $C(L)\neq C(K)$. Suppose $P\in C(L)$. I say that $P$ is an $L$-\emph{primitive point} if $L/K$ is the smallest Galois extension of $K$ over which $P$ is defined. 

Note that all the examples constructed in this subsections are of curves equipped with primitive points:  the condition of Theorem~\ref{th:main0}\ref{th-points-3} says that $P=(u,\pm d)$ is an $L_f$-primitive point.


I make the following very optimistic conjecture (and even more optimistic if one wants to remove the hyperelliptic condition):
\begin{conj}\label{con:main}
Let $C/K$ be a smooth, projective, hyperelliptic curve of genus $g\geq 1$.  Suppose $C$ has a $K$-rational point $P_0$ and $C$ is embedded in $J_C$ using this $K$-rational point. Suppose $L/K$ is a finite Galois extension  with $G=\gal(L/K)$ and $P\in C(L)$ is an $L$-primitive point. If $G$ is a nonabelian almost simple group and $|G|$ is sufficiently large then $P$ is non-torsion in $J_C(L)$; in particular $$\dim J_C(L)\tensor\Q\geq  
m(\gal(L/K)).$$ 
\end{conj}

A particular case of the conjecture is the following conjecture.
\begin{conj}\label{con:main-osada-polynomial}
Let $n\geq 5$ be an odd integer. Let $0\neq d\in\Z$ be an integer chosen such that $g(X)=X^n-X-1+d^2\in\Q[X]$ has distinct roots. Let
$C/\Q$ be a smooth, projective, hyperelliptic curve of genus $g\geq 2$ defined by $Y^2=g(X)$ and  embed in $C\into J_C$ using the tautological point at $\infty$.
Let $L_f/\Q$ be the splitting field of $f(X)=X^n-X-1$. 

Then for any root $u$ of $f(X)=0$, the $L_f$-rational point  is non-torsion in $J_C(L)$ and in particular $$\dim J_C(L_f)\tensor\Q\geq  
	m(\gal(L_f/\Q))=n-1.$$ 
\end{conj}

\subsection{Two numerical examples of the Conjecture}\label{ss:examples}
Here are two concrete example of  Conjecture~\ref{con:main}: 
\begin{example}
Consider $f(X)=x^5-x-1\in\Q[X]$ and $g(X)=f(x)+7^2$ (so $d=7$). By \cite{osada87} the Galois group of $f(X)$ (over $\Q$) is $S_5$. The curve $C:Y^2=g(X)=X^5-X-1+7^2$ has genus $g=2$ and the curve has rational points $(u,\pm 7)$ over the splitting field $L_f$ of $f(X)$ where $u$ is any root of $X^5-X-1=0$. To check if $(u,7)$ is non-torsion it is sufficient to work in the field $\Q(u)$ obtained by attaching one root $u$ of $f(X)$ (this is of degree five). One can calculate (using \cite{sage}) multiples of this point in the Jacobian of $C$ and as far as I see the points $(u,7)\in C(L_f)$ are non-torsion in the Jacobian of $C$. So $J_C(L_f)\tensor\Q$ is equipped with a faithful, $\Q$-rational representation of $S_5$ and hence has $\Q$-rank at least four.  
\end{example}

\begin{example}
Here is a genus three example with a non-abelian simple group. It is well-known that the polynomial $f(x)=x^7-7x+3$ has as its Galois group over $\Q$ the unique simple group of order $168$ (isomorphic to $PSL(2,7)$). Consider the curve $C:y^2=x^7-7x+3+11^2$ so $d=11$ is chosen such that the polynomial on the right has distinct roots. Let $L/\Q$ be the splitting field of $x^7-7x+3$. Then $\gal(L/\Q)=PSL(2,7)$ is simple of order $168$. By construction $C(L)$ has at least $14$ points given by $(u,\pm d)$ where $u$ is any root of $x^7-7x+3$. One expects that the rank of $J_C(L)$ is at least $6$ as this is the dimension of the smallest faithful $\Q$-rational representation of $\gal(L/\Q)=PSL(2,7)$. By computing with \cite{sage} one checks that by adjoining any root $u$ of $f(x)$ to $\Q$, the point $(u,11)$ is not of order $\leq 100$ over $\Q(u)$ and it is quite reasonable to expect that this point is not torsion and hence $J_C(L)\tensor Q$ has $\Q$-rank at least six.
\end{example}

\subsection{Proof of the rank conjecture in some cases}
Let me now prove Conjecture~\ref{con:main} in some cases.
\bthm\label{th:main1}
Let $f(X)\in\Z[X]$ be a monic irreducible polynomial of odd degree $n\geq 5$ and with Galois group $G=\gal(L_f/\Q)$ which is either simple or  $G=S_n$. Let $L_f$ be a  splitting field of $f(X)$. Suppose $\wp$ is a non-zero prime ideal in $L_f$ not lying over $(2)\subset \Z$. Suppose that $d\in L_f$ is such that
\benum[label={\bf(\arabic{*})}]
\item $d\in\wp$,
\item $g(X)=f(X)+d^2\in L_f[X]$ is irreducible, and
\item $g(X)\cong f(X)\bmod{\wp}$ has distinct roots,
\item let $C: Y^2=g(X)=f(X)+d^2$ be the hyperelliptic curve over $L_f$ defined by $g(X)$,
\item embed $C\into J=J_C$ using $\infty\in C(L_f)$.
\eenum
Then 
\benum[label={\bf(\arabic{*})}]
\item for any root $u$ of $f(X)=0$, $C$ has $L_f$-rational points $(u,\pm d)$, 
\item for each $u$, these points are of infinite order in the Jacobian $J$ of $C$,
\item and $J(L_f)$ has rank at least $m(G)$ (resp. $n-1=2g$ if $G=S_n$).
\eenum
\ethm
\bp 
The proof is by methods of \cite{joshi99}. By Theorem~\ref{th:main0} it suffices to prove that the listed points are not torsion in the Jacobian. As $d\in\wp$ and $g(X)\cong f(X)\bmod{\wp}$ has distinct roots, the images of these points modulo $\wp$ generate a subgroup of the $2$-torsion of this curve modulo $\wp$ (and so the torsion order is prime to the residue characteristic of $\wp$). In particular the $n$ points $(u,d)$  are of order at most two in the Jacobian. But $J[2]$ is generated by the $2g+1=n$ Weierstrass points of $C$ which correspond to the roots of the polynomial $g(X)$. As $g(X)$ is irreducible (over $L_f$),  no two torsion point of $J$ is defined over $L_f$. So the listed points are of infinite order in $J(L_f)$. If $G$ is simple then by Theorem~\ref{th:main0} the rank is at least $m(G)$. So now assume $G=S_n$. As the action of $G=\gal(L_f/\Q)$ on the points $n\geq 5$ points $(u,d)$ is transitive,  this action cannot factor through $S_n/A_n=\Z/2$. So the representation of $G$ on $J(L_f)\tensor\Q\neq 0$ is faithful and as $m(S_n)=n-1$ for $n\geq 5$, the rank of $J(L_f)$ is at least $\geq  n-1=2g$ if $G=S_n$ as was claimed. 
\ep

\subsection{Elliptic Curve variant with points over simplest cubic fields}\label{ss:cubic}
My next example, a variant of both \eqref{eq:orig} and the equation in Theorem~\ref{th:main0}, is inspired by \cite{shanks74} and provides a simple example of the abelian group case of the preceding examples. One can combine methods of \cite{shanks74} with those outlined above to produce elliptic curves over $\Q$ with three rational points defined over the ``simplest cubic fields'' constructed in \cite{shanks74}. Briefly the construction is as follows. Suppose $a$ is a positive integer and $$f(X)=X^3-(aX^2+(a+3)X+1).$$ Let $L_f$ be the splitting field of $f(X)\in\Z[X]$. In \cite{shanks74} it is shown that $\gal(L_f/\Q)\isom A_3=\Z/3$ and these fields are called simplest cubic fields. If $\rho$ is any root of $f(X)$ then other roots of $f(X)$ are $\rho_2=-1/(1+\rho), \rho_3=-1/(1+\rho_2)$. Moreover $\rho,\rho_2,\rho_3$ are fundamental units in $L_f$. 

Now consider the elliptic curves $E: Y^2=f(X)+d^2$. Then this has $L_f$-rational points $(\rho,\pm d),(\rho_2,\pm d)$ and $(\rho_3,\pm d)$. 

Here two numerical examples.

Let $a=50$ then $L_f/\Q$ has class number $h=19$. Choose $d=50$ then the torsion subgroup of $E(L_f)$ is trivial (as can be computed using \cite{sage}) so these points are necessarily of infinite order and one finds the height pairing matrix for the two $L_f$-rational points $P=(u,50)$ and $Q=(-1/(1+u),50)$:
\be 
\begin{pmatrix}
	3.49761746081572 & -1.74880873040786\\
	-1.74880873040786 & 3.49761746081572
\end{pmatrix}
\ee
whose determinant is approximately $9.17499592665226$. 

Now set $R=(\rho_3,50)$ and $P_1=P+R,P_2=P-R$. Then the height pairing matrix of $P_1,P_2$ is
\be 
\begin{pmatrix}
	3.49761746081572 & 0.000000000000000\\
	0.000000000000000 & 10.4928523824472
\end{pmatrix}
\ee
So the rank of $E(L_f)$ is at least two and $P_1,P_2$ are orthogonal with respect to the height pairing.

Here is the second example.

Now let $a=36$ then $L_f/\Q$ has class number $h=12$ (this case is not included in the tables in \cite{shanks74}). Now choose $d=a=36$ then the torsion subgroup of $E(L_f)$ is trivial (as can be computed using \cite{sage}) so these points are necessarily of infinite order and one finds the height pairing matrix for the two $L_f$-rational points $P_1=P+Q$ and $P_2=P-Q$ where $P=(u,36)$ and $Q=(-1/(1+u),36)$ 
\be 
\begin{pmatrix}
  2.49439281440290 &  8.88178419700125\times 10^{-16}\\
8.88178419700125\times 10^{-16} &     7.48317844320869
\end{pmatrix}
\ee
So the rank of $E(L_f)$ is at least two and $P_1,P_2$ are orthogonal with respect to the height pairing.

Hence it seems reasonable to conjecture that for any positive integers $a>2$ such that $g(X)=f(X)+a^2$ has distinct roots the rank of $E(L_f)$ is at least two and the  points $P_1,P_2$ as constructed above provide an orthogonal basis for the $\Q$-vector space spanned by the points $P,Q,R$ in $E(L_f)$.

In fact it seems reasonable to expect that for any positive integers $a,d$ such that $f(X)+d^2$ has distinct roots, the curve $E$ has rank at least two over the cubic (galois) extension $L_f/\Q$. Then $P_1=P+Q,P_2=P-Q$ as above are orthogonal for the height pairing. Let me also note that the hypothesis $a,d\in\Z$ cannot be relaxed as the following example shows.

Let $a=745, d=u^2 - u - 3$ (both chosen randomly by SAGE) then the height pairing matrix for $P_1,P_2$ is
$$
\begin{pmatrix}
6.40170943412354 & -1.10204026919419\\
-1.10204026919419 & 17.0158570702399
\end{pmatrix} 
$$
So orthogonality of $P_1,P_2$  is certainly false if one takes $d\in L_f$ (in fact $d=u^2-u-3$ is an algebraic integer). Moreover $P,Q$ are not orthogonal either and orthogonality also fails if one chooses $a\in\Q$ but not in $\Z$ and hence orthogonality seems to require $a\in\Z$ and $d\in\Q$.

Here is what I can prove using methods of \cite{joshi99}.

\bthm\label{th:main2}
Let $a\geq 1$ be an integer and let $f(X)=X^3-(aX^2+(a+3)X+1)$ with roots $\rho,\rho_2,\rho_3$. Let $L_f$ be its splitting field over $\Q$. Suppose $\wp$ is a non-zero prime ideal in $L_f$ not lying over $(2)\subset \Z$. Suppose that $d\in L_f$ is such that
\benum[label={\bf(\arabic{*})}]
\item $d\in\wp$,
\item $g(X)=f(X)+d^2\in L_f[X]$ is irreducible, and
\item $g(X)\cong f(X)\bmod{\wp}$ has distinct roots.
\eenum
Then 
\benum[label={\bf(\arabic{*})}]
\item $E$ has $L_f$-rational points $P_1=(\rho,\pm d),P_2=(\rho_2,\pm d), P_3=(\rho_3,\pm d)$, 
\item each of these points is of infinite order,
\item and $E(L_f)$ has rank at least two.
\eenum
\ethm
\bp 
The proof is similar to that of Theorem~\ref{th:main1} and is left to the reader. 
\ep

\begin{conj}\label{con:orthogonality}
Let $a\geq 1$ be any integer and $0\neq d\in\Q$ such that $f(X)+d^2$ has distinct roots. In the notation of  Theorem~\ref{th:main2} and with $P_1=(\rho,d,1), P_2=(\rho_2,d,1), P_3=(\rho_3,d,1)$ and $P_i+P_j$, $P_i-P_j$, (for $1\leq i\neq j\leq 3$) are always orthogonal for the canonical height pairing on the Mordell-Weil  group $E(L_f)$. In particular $E(L_f)$ has rank at least two.
\end{conj}

Let me remark that orthogonality is immediate if one proves that $P_1,P_2,P_3$ span a vector space of dimension two (note that $P_1,P_2,P_3$ satisfy one linear relation $P_1+P_2+P_3=O$) upon observing that $P_1,P_2,P_3$ have the same N\'eron-Tate height (as they are conjugate).
%
\bibliographystyle{plain}
\bibliography{points.bib,../../master/joshi.bib,../../master/master6.bib}
\end{document}